\documentclass{gtart_a}
\pdfoutput=1

\usepackage{pinlabel}
\usepackage[all]{xy}

\title{Duality for Legendrian contact homology}

\author{Joshua M Sabloff} 
\givenname{Joshua M}
\surname{Sabloff}
\address{Haverford College\\Haverford, PA 19041\\USA} 
\email{jsabloff@haverford.edu}
\urladdr{}

\volumenumber{10}
\issuenumber{}
\publicationyear{2006}
\papernumber{53}
\lognumber{0652}
\startpage{2351}
\endpage{2381}

\doi{}
\MR{}
\Zbl{}


\received{7 September 2005}
\revised{29 July 2006}
\accepted{9 November 2006}
\proposed{Yasha Eliashberg}
\seconded{Leonid Polterovich, Joan Birman}
\published{8 December 2006}
\publishedonline{8 December 2006}
\corresponding{}
\editor{CPR}
\version{}

\arxivreference{math.SG/0508187}

\let\xysavmatrix\xymatrix
\def\xymatrix{\disablesubscriptcorrection\xysavmatrix}
\AtBeginDocument{\let\bar\wbar\let\tilde\wtilde\let\hat\what}

\newcommand{\df}{\ensuremath{\partial}}
\newcommand{\dfe}{\ensuremath{\partial^\varepsilon}}

\newcommand{\fdf}{\ensuremath{\widehat{\partial}}}

\newcommand{\alg}{\ensuremath{\mathcal{A}}}
\newcommand{\falg}{\ensuremath{\widehat{\mathcal{A}}}}

\DeclareMathOperator{\eval}{ev}
\newcommand{\fc}{\ensuremath{\kappa}} 

\newcommand{\rr}{\ensuremath{\mathbb{R}}}
\newcommand{\zz}{\ensuremath{\mathbb{Z}}}

\makeatletter
\def\cnewtheorem#1[#2]#3{\newtheorem{#1}{#3}[section]
\expandafter\let\csname c@#1\endcsname\c@thm}

\theoremstyle{plain}
\newtheorem{thm}{Theorem}[section]
\cnewtheorem{cor}[thm]{Corollary}
\cnewtheorem{lem}[thm]{Lemma}
\cnewtheorem{claim}[thm]{Claim}
\cnewtheorem{prop}[thm]{Proposition}
\cnewtheorem{conj}[thm]{Conjecture}

\theoremstyle{definition}
\cnewtheorem{defn}[thm]{Definition}

\theoremstyle{remark}
\cnewtheorem{rem}[thm]{Remark}
\cnewtheorem{exa}[thm]{Example}

\numberwithin{equation}{section}

\makeatother  
\makeautorefname{defn}{Definition}

\begin{document}

\begin{asciiabstract}
The main result of this paper is that, off of a `fundamental class' in
degree 1, the linearized Legendrian contact homology obeys a version
of Poincare duality between homology groups in degrees k and -k.  Not
only does the result itself simplify calculations, but its proof also
establishes a framework for analyzing cohomology operations on the
linearized Legendrian contact homology.
\end{asciiabstract}

\begin{abstract}
  The main result of this paper is that, off of a ``fundamental
  class'' in degree 1, the linearized Legendrian contact homology
  obeys a version of Poincar\'e duality between homology groups in
  degrees $k$ and $-k$.  Not only does the result itself simplify
  calculations, but its proof also establishes a framework for
  analyzing cohomology operations on the linearized Legendrian contact
  homology.
\end{abstract}

\maketitle

\section{Introduction}
\label{sec:intro}

\subsection{Legendrian contact homology}
\label{ssec:lch}

As in smooth knot theory, a fundamental problem in Legendrian knot
theory is to find effective invariants and to understand their
structure and meaning.  Bennequin \cite{bennequin} initiated the
modern study of Legendrian knots by introducing two ``classical''
invariants: the \textit{Thurston--Bennequin number} (which measures
the difference between the framing of a knot coming from the contact planes and
the Seifert surface framing) and the \textit{rotation number} (which
measures the twisting of the tangent to the knot inside the contact
planes with respect to a suitable trivialization).  These two
invariants suffice to classify Legendrian knots in the standard
contact structure on $\rr^3$ when the underlying smooth knot type is
the unknot (Eliashberg and Fraser \cite{yasha-fraser}), a torus knot
or the figure eight knot (Etnyre and Honda \cite{etnyre-honda:knots})
or a cable link (Ding and Geiges \cite{ding-geiges:knots}).

The first non-classical invariant of Legendrian knots was
\textit{Legendrian contact homology}, a Floer-type theory that comes
from geometric ideas of Eliashberg and Hofer \mbox{\cite{egh,yasha:icm}} and
was rendered combinatorially computable by Chekanov \cite{chv} for
knots in the standard contact $\rr^3$. The Legendrian contact homology
of a Legendrian knot is the homology of a freely-generated
differential graded algebra (DGA) $(\alg, \df)$, which we shall refer
to as the \textit{Chekanov--Eliashberg DGA}, itself an invariant up to
``stable tame isomorphism.'' It is difficult to extract information
from the stable tame isomorphism class of the Chekanov--Eliashberg
DGA, but Chekanov defined a linearized version which was sufficient to
distinguish the first examples of Legendrian knots with the same
smooth knot type and classical invariants \cite{chv}.  The homology of
the linearized DGA is usually encoded in a
\textit{Poincar\'e--Chekanov polynomial}, in which the coefficient of
$t^k$ denotes the dimension of the linearized homology in degree $k$.
The set of Poincar\'e--Chekanov polynomials taken over all possible
linearizations of the DGA is invariant under Legendrian isotopy.
Though other more powerful --- but less computable --- methods of
distilling information from the Chekanov--Eliashberg DGA have been
developed (Ng \cite{lenny:computable}), there is still much to be
learned about Chekanov's linearized theory.

Many recent advances in Floer-type theories, and contact homology in
particular, have come from importing classical Morse-theoretic ideas.
For example, K\'alm\'an's analysis of nontrivial loops of Legendrian
knots was motivated by continuation maps in Morse theory
\cite{kalman:mono1}.  Another example is the extension of the
combinatorial definition of the Chekanov--Eliashberg DGA to Legendrian
knots in circle bundles with contact structures transverse to the
fiber, which was achieved by transporting Morse--Bott methods into the
contact homology picture (Sabloff \cite{s1bundles}); see also
Bourgeois' work on Morse--Bott methods for the non-relative version of
contact homology \cite{bourgeois:mb}.

This paper translates ideas of Betz and Cohen \cite{cohen-betz} and
Fukaya and Oh \cite{fukaya-oh} on gradient flow trees into the
construction of a Poincar\'e duality map, a cap product, and a
fundamental class for the linearized DGA.  The main result is the
following duality theorem:

\begin{thm} \label{thm:pd} If $h_n$ denotes the dimension in degree
  $n$ of the homology of a linearization of the Chekanov--Eliashberg
  DGA of a Legendrian knot, then:
  \begin{align*}
    h_n &= h_{-n} & |n|>1, \\
    h_1 &= h_{-1}+1. &
  \end{align*} 
\end{thm}

Said another way, the Poincar\'e--Chekanov polynomial $P(t)$ of a
linearization satisfies:
\begin{equation*}
  P(t) = P(t^{-1}) + (t - t^{-1}).
\end{equation*}
The extra class in degree $1$ will turn out to be the fundamental
class, which is closely related to the fundamental class of the
circle. This theorem can greatly simplify calculations (see Melvin and
Shrestha \cite{melvin-shrestha}, for example), but of greater interest
is its proof, which introduces an algebraic and geometric framework
for understanding the ``algebraic topology'' of the linearized DGA.

\subsection{Morse-theoretic motivation}
\label{ssec:morse}

While the proof of \fullref{thm:pd} will be combinatorial in
nature, it is motivated by geometric ideas from Morse theory.
Classically, the Poincar\'e duality map caps a cohomology class with
the fundamental class.  Using gradient flow trees as in
\cite{cohen-betz,fukaya-oh}, the cap product and fundamental class
can be reinterpreted in the setting of Morse--Witten theory.  Recall
that the Morse--Witten complex $C_*(M,f)$ of a manifold $M$ with Morse
function $f$ is generated by the critical points of $f$, while the
differential comes from counting rigid negative gradient flow lines
between critical points; see Schwarz \cite{schwarz}. The cochain complex
$C^*(M,f)$ is generated by the same critical points, but the
codifferential counts \emph{upward} gradient flows.

\begin{figure}[ht!]
  \centerline{\includegraphics[scale=.9]{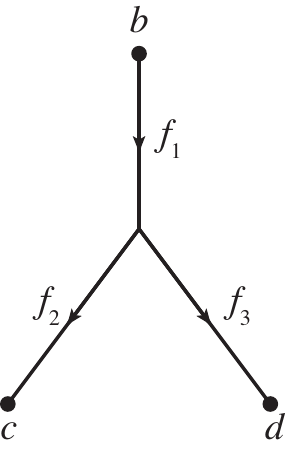}}
  \caption{A gradient flow tree that gives the cap product $[d] =
    \gamma \cap B$ between $B=[b]$ and $\gamma=[c]$.  Notice that
    there is one homological input at $b$, one cohomological input at
    $c$, and a homological output at $d$.}  \label{fig:Y}
\end{figure}

In this language, the fundamental class is represented by a sum of the
maxima of $f$.  More interestingly, the cap product between a homology
class $B$ and a cohomology class $\gamma$ is computed by counting
certain rigid ``Y''--shaped gradient flow trees.  Choose three Morse
functions $f_1$, $f_2$, and $f_3$, and choose representatives $b$ for
$B$ in $C_*(M,f_1)$ and $c$ for $\gamma$ in $C^*(M, f_2)$.  The
``Y''--shaped tree follows the negative gradient flow of $f_1$ out of
$b$ and splits at the vertex into negative gradient flows for $f_2$
and $f_3$.  These two flowlines end at $c$ and at a critical point of
$f_3$ that represents $\gamma \cap B$ in $C_*(M,f_3)$, respectively;
see \fullref{fig:Y}.  More generally, gradient flow trees give rise
to cohomology operations, with homology inputs and cohomology outputs
at the outward-flowing, or positive, ends of a negative gradient flow
line, and homology outputs and cohomology inputs at the negative ends.
Of particular interest in this paper is the Poincar\'e duality
isomorphism, which, as Betz and Cohen note, comes from a tree with two
positive ends and one vertex, where one of the positive ends is a
homology input and the other a cohomology output.

Legendrian contact homology fits into this picture as the
Morse--Witten--Floer theory of the action functional on the space of
paths that begin and end on a Legendrian knot $K$.  In this setting,
Reeb chords are the critical points that generate the
Morse--Witten--Floer complex, while holomorphic disks in the
symplectization $(\rr^3 \times \rr, d(e^t \alpha))$ take the place of
the gradient flow lines that define the differential \df.  In order to
ensure that $\df^2=0$, the holomorphic disks must be allowed to have
multiple negative ends.  The end result is a non-commutative DGA
generated by the Reeb chords.  Chekanov's combinatorial formulation of
the DGA arises from the correspondence between Reeb chords and
crossings in the $xy$ projection of a Legendrian knot --- the Reeb
direction in the standard contact $\rr^3$ is the $z$ direction --- and
between holomorphic disks in the symplectization and immersed disks in
the $xy$ projection; see \fullref{sec:background} for the
combinatorics and Etnyre, Ng and Sabloff \cite{ens} for more details of
the translation between geometry and combinatorics.

The Morse-theoretic basis for Legendrian contact homology would lead
one to expect that the longest Reeb chords should give a fundamental
class in Legendrian contact homology, and the cap product should come
from disks with one positive and two negative ends.  The analogy
between Legendrian contact homology and Morse theory breaks down in
the case of the Poincar\'e duality map: disks with multiple positive
ends do not appear in the contact homology theory.  To access disks
with two positive ends, it is necessary to expand the algebraic
framework of the Chekanov--Eliashberg DGA.  The natural expansion is a
relative version of Eliashberg, Givental, and Hofer's symplectic field
theory (SFT) \cite{egh}.  The definition of a ``Legendrian SFT,''
however, runs into some subtle issues regarding compactness of the
moduli space of curves used to define the differential and has yet to
be rigorously defined. Instead, the expansion can be achieved by an
appropriate interpretation of the generators and differentials of the
DGA of a link formed by several vertical translates of the original
knot.

The remainder of the paper is organized as follows: background notions
--- including descriptions of Legendrian knot diagrams, the
Chekanov--Eliashberg DGA, and the linearization procedure --- are
described in \fullref{sec:background}. Next,
\fullref{sec:expanded-alg} describes the structure of the
``expanded DGA.''  \fullref{sec:duality} describes the
relationship between the original linearized DGA and the linearized
expanded DGA, resulting in the definition of the Poincar\'e duality
map and the proof of \fullref{thm:pd}.  Finally,
\fullref{sec:alg-top} explores the ``algebraic topology'' of the
linearized DGA, defining both the fundamental class and the cap
product and showing that a ``capping with the fundamental class''
construction inverts the duality isomorphism.

\subsubsection*{Acknowledgments}
This paper greatly benefited from discussions with John Etnyre, Paul
Melvin, Lenny Ng, and Lisa Traynor. The referee's insightful
comments spurred me to clean up some ideas surrounding the fundamental
class, which, in turn, led to a greatly improved proof of duality. The
research behind this paper started out in a rather different
direction, and Yasha Eliashberg and Fr\'ed\'eric Bourgeois were
instrumental in helping me understand the ideas of symplectic field
theory.

\section{Background notions}
\label{sec:background}

\subsection{Legendrian knots and their diagrams}
\label{ssec:knots}

This section briefly reviews the basic notions of Legendrian knot
theory; for a more comprehensive introduction, see
Etnyre \cite{etnyre:knot-intro} or Sabloff \cite{lecnotes}.

The \textit{standard contact structure} on $\rr^3$ is the completely
non-integrable $2$--plane field given by the kernel of $\alpha = dz -
y\,dx$.  A \textit{Legendrian knot} is a smooth embedding $K\co  S^1 \to
\rr^3$ that is everywhere tangent to the contact planes.  That is, the
embedding satisfies $\alpha(K') = 0$.  An ambient isotopy of $K$
through other Legendrian knots is a \textit{Legendrian isotopy}.
Legendrian knots are plentiful; for example, any smooth knot can be
continuously approximated by a Legendrian knot.

There are two useful projections of Legendrian knots.  The first is
the \textit{front projection} $\pi_f$ to the $xz$ plane. In the front
projection, the $y$ coordinate of a knot may be recovered from the
slope of its projection.  As a result, the projection can have no self
or vertical tangencies; it has semi-cubical cusps instead.  Further,
the crossing information is completely determined: the strand with
lesser slope will always pass in front of the strand with greater
slope. Any circle in the $xz$ plane that has no vertical tangencies,
has no self-tangencies, and that is immersed except at finitely many
cusps lifts to a Legendrian knot.

\begin{figure}[ht!]\small
  \centerline{\includegraphics{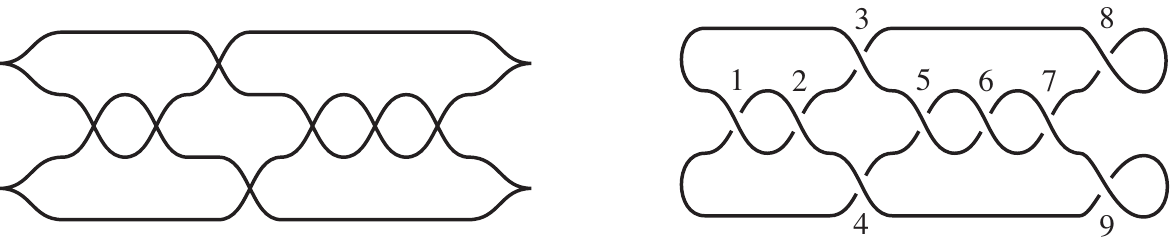}}
  \cl{(a)\hskip 2in(b)} 
\caption{The (a) front and (b) Lagrangian diagrams of a Legendrian $5_2$
    knot. The meaning of the numbers will become clear in
    \fullref{ssec:ch-r3}.}  
   \label{fig:numbered-5-2}
\end{figure}

Though the front projection is easier to use, it is more natural to
define Legendrian contact homology using the \textit{Lagrangian
  projection} $\pi_l$ to the $xy$ plane.  Unlike the front projection,
not every immersion into the $xy$ plane is the Lagrangian projection
of a Legendrian knot $K$: a system of inequalities involving the areas
of the connected components $\rr^2 \setminus \pi_l(K)$ must be
satisfied (see \cite{chv}).  It is simpler (and sufficient) to work
with \textit{Lagrangian diagrams} of $K$, ie, immersions $D$ of the
circle into the $xy$ plane, together with crossing information, for
which there is an orientation-preserving diffeomorphism of the plane
carrying $D$ to $\pi_l(K)$.  See \fullref{fig:numbered-5-2}(b), for
example.

Ng's \textit{resolution} procedure (see \cite{lenny:computable}) gives
a canonical translation from a front projection to a Lagrangian
diagram.  Combinatorially, there are three steps:
\begin{enumerate}
\item smooth the left cusps;
\item replace the right cusps with loops (see the right side of the
  Lagrangian projection in \fullref{fig:numbered-5-2}); and
\item resolve the crossings so that the overcrossing is the one with
  lesser slope.
\end{enumerate}
A key feature of the resolution procedure is that the heights of the
crossings in the Lagrangian diagram strictly increase from left to
right, with the jumps in height between crossings as large as desired.
In particular, the crossings in a resolved Lagrangian diagram that
come from the rightmost cusps have the greatest height among all
crossings.

\subsection{The Chekanov--Eliashberg DGA}
\label{ssec:ch-r3}

The DGA was originally defined by Chekanov in \cite{chv} for
Lagrangian diagrams; see also \cite{ens,lenny:computable,lecnotes}.
This section contains a brief review of the definition.

Let $K$ be an oriented Legendrian knot in the standard contact $\rr^3$
with a generic Lagrangian diagram $D$.  Label the crossings with
$\{q_1, \ldots, q_n\}$.  Let $A$ be the graded vector space over
$\zz/2$ generated by $q_1, \ldots, q_n$, and let $\alg$ be the graded
free unital tensor algebra $TA$.

The grading is determined by the assignment of a \textit{capping path}
to each crossing.  A capping path is one of the two immersed paths
that start at the overcrossing of $q_i$, trace out a portion of $D$,
and end when $D$ first returns to $q_i$, necessarily at an
undercrossing.  Assume, without loss of generality, that the strands
of $D$ at each crossing are orthogonal. The \textit{grading} of $q_i$
is defined to be:
\begin{equation*} \label{eqn:grading}
  |q_i| \equiv 2 r(\gamma_i) - \frac{1}{2} \mod 2r(K).
\end{equation*}
Extend the grading to all words in $\alg$ by letting the grading of a
word be the sum of the gradings of its constituent generators.

\begin{figure}[ht!]
  \centerline{\includegraphics[scale=.9]{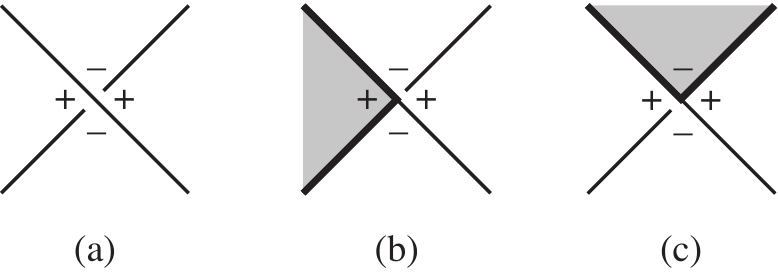}}
  \caption{(a) A labeling of the quadrants surrounding a crossing; (b)    
    The image of $f \in \Delta$ near the crossing $q_i$ and (c) the
    crossings $q_{j_l}$. The other positive (resp.\ negative) corner in
    (b) (resp.\ (c)) is also possible.}  \label{fig:disks}
\end{figure}

The next step is to define a differential on $\alg$ by counting
immersions of the $2$--disk $B^2$ with an ordered set of marked points
on its boundary. The number of marked points must be at least one, but
otherwise can vary from disk to disk; if there are $k+1$ of marked
points, label them $\{z_0, \ldots, z_k\}$.  Further, label the corners
of $D$ as in \fullref{fig:disks}(a).  The immersions of interest
are the following:

\begin{defn}
  \label{defn:ce-immersed} Given an ordered set of generators $q_i,
  q_{j_1}, \ldots, q_{j_k}$, define the set $\Delta(q_i; q_{j_1},
  \ldots, q_{j_k})$ to consist of orientation-preserving immersions
  \begin{equation*}
    f\co  B^2 \to \rr^2
  \end{equation*}
  up to smooth reparametrization that map $\partial B^2$ to the image
  of $D$ subject to the following conditions:
  \begin{enumerate}
  \item The restriction of $f$ to the boundary is an immersion away
    from the marked points $z_i$.
  \item The map $f$ has the property that $f(z_0) = q_i$ and $f(z_l) =
    q_{j_l}$, such that the points $q_i, q_{j_1}, \ldots, q_{j_k}$ are
    encountered in counter-clockwise order along the boundary.
  \item In a neighborhood of the points $q_i, q_{j_1}, \ldots,
    q_{j_k}$, the image of the disk under $f$ has the form indicated
    in \fullref{fig:disks}(b) and (c).
\end{enumerate}
\end{defn}

Note the analogy between the positive and negative corners in
\fullref{fig:disks} and positive and negative ends of a gradient
flow tree.  With this definition in hand, define the differential 
as follows:

\begin{defn}\label{defn:ce-diffl}\qquad$\displaystyle
\df(q_i) = \sum_{\text{Words } q_{j_1} \cdots q_{j_k}} 
\sum_{f \in \Delta(q_i; q_{j_1}, \ldots,
q_{j_k})} q_{j_1} \cdots q_{j_k}$

  Extend \df\ to all of \alg\ via linearity and the Leibniz rule.
\end{defn}

The fact that the sum in the definition of $\df$ is finite comes from
the following lemma, which is essentially an application of Stokes'
Theorem:

\begin{lem} \label{lem:stokes}
  Let $h(x)$ be the height of the crossing $x$.  If there is an
  immersed disk that satisfies all but the third condition of
  \fullref{defn:ce-immersed} and has positive corners at $x_1,
  \ldots, x_k$ and negative corners at $y_1, \ldots, y_l$, then
  \begin{equation*}
    \sum_i h(x_i) - \sum_i h(y_i) > 0.
  \end{equation*}
\end{lem}

The central results in the theory are:

\begin{thm}[Chekanov \cite{chv}] \
  \label{thm:dga}
  \begin{enumerate}
  \item The differential \df\ has degree $-1$.
  \item The differential satisfies $\df^2=0$.
  \item The ``stable tame isomorphism class'' (and hence the homology)
    of the DGA is invariant under Legendrian isotopy.
  \end{enumerate}
\end{thm}

Here, a \textit{stabilization} of \alg\ is a DGA $S(\alg)$ with two
new generators $a$ and $b$ such that $\df a = b$ and $\df b = 0$.  Two
DGAs \alg\ and $\alg'$ are \textit{stable tame isomorphic} if there
exist (possibly multiple) stabilizations of each that are tame
isomorphic. See \cite{chv} for a description of the technical
condition ``tame;'' this condition allows for the invariance of
\emph{based} algebras and is essential for the proof of the invariance
of the linearized contact homology described in the next section,
though its exact technical definition is not necessary for the
arguments in this paper.

\begin{exa} \label{ex:5-2-diff}
  The knot in \fullref{fig:numbered-5-2} has nine crossings.  The
  generators $q_1$ and $q_2$ have grading $0$; $q_3$, $q_4$, $q_8$,
  and $q_9$ have grading $1$; $q_5$ and $q_7$ have grading $2$; and
  $q_6$ has grading $-2$.  The differential is given by:
  \begin{equation}
    \begin{split}
      \df q_1 &= \df q_2 = \df q_6 = \df q_7 = 0 \\
      \df q_3 &= 1+ q_1 q_2 \\
      \df q_4 &= 1+ q_2 q_1 \\
      \df q_5 &= q_3 q_1 + q_1 q_4 \\
      \df q_8 &= 1+ q_1 + q_1 q_6 q_7 \\
      \df q_9 &= 1+ q_1 + q_7 q_6 q_1.
    \end{split}
  \end{equation}
\end{exa}

\subsection{Linearized contact homology}
\label{ssec:linear}

Chekanov introduced an important computable invariant of the stable
tame isomorphism class of the Chekanov--Eliashberg DGA called the
\textit{linearized contact homology}.  The differential \df\ on \alg\
may be split into a sum of differentials $\df = \sum_{l=0}^\infty
\df_l$, where the image of $\df_0$ lies in the ground ring $\zz/2$ and
$\df_l$ maps a generator of $A$ into $A^{\otimes l}$.  If $\df_0 = 0$,
then the equation $\df^2 = 0$ implies $\df_1^2 = 0$.  In this case, it
follows that $(A, \df_1)$ is an honest chain complex whose homology
$H_*(A, \df_1)$ may be easily computed.

Rarely does the DGA have the property that $\df_0 = 0$.  Suppose,
however, that there exists an algebra map $\varepsilon\co \alg \to
\zz/2$ that satisfies:
\begin{enumerate}
  \item $\varepsilon \circ \df = 0$, and
  \item $\varepsilon (q_i) = 0$ if $|q_i| \neq 0$.
\end{enumerate}
Such a map is called an \textit{augmentation}; if $ \varepsilon(q_i) =
1$, say that $q_i$ is \emph{augmented}.  Augmentations are not
uncommon, but they do not always exist; it turns out that their
existence is equivalent to the existence of a ``ruling'' of a front
diagram of a Legendrian knot (Fuchs \cite{fuchs:augmentations}, Fuchs
and Ishkhanov \cite{fuchs-ishk} and Sabloff \cite{rulings}).

Given an augmentation $\varepsilon$, define an automorphism
$\Phi^\varepsilon\co  \alg \to \alg$ by:
\begin{equation*}
  \Phi^\varepsilon(q_i) = q_i + \varepsilon(q_i).
\end{equation*}
Let \dfe\ be the differential induced by $\Phi^\varepsilon$:  
\begin{equation*}
  \dfe = \Phi^\varepsilon  \df  (\Phi^\varepsilon)^{-1}.
\end{equation*}
It is straightforward to verify that the DGA $(\alg, \dfe)$ satisfies
$\dfe_0 = 0$.

\begin{rem}
  Once the augmentation is known, it is possible to read off the
  linearized differential directly from the Lagrangian diagram.  An
  immersed disk that contributes $q_j$ to the linearized differential
  of $q_i$ has a positive corner at $q_i$, a negative corner at $q_j$,
  and possibly other negative corners at augmented crossings. For
  example, if $\df q_1 = q_2 q_3$ and only $q_2$ is augmented, then
  $\dfe_1 q_1 = q_3$.  If both $q_2$ and $q_3$ are augmented, then
  $\dfe_1 q_1 = q_2 + q_3$.\end{rem}

For each augmentation, there is a \textit{Poincar\'e--Chekanov
  Polynomial}:
\begin{equation*}
  P_{\varepsilon}(t) = \sum_{n=-\infty}^\infty \dim H_n(A,
  \dfe_1) \cdot t^n.
\end{equation*}
Chekanov proved in \cite{chv} that the set $\mathcal{P} = \left\{
  P_{\varepsilon}(t) \right\}_{\varepsilon \in \mathcal{E}}$ is
invariant under Legendrian isotopy, where $\mathcal{E}$ is the set of
all possible augmentations of $(\alg, \df)$. The proof of invariance
splits into two parts: the first is the somewhat trivial verification
that $\mathcal{P}$ does not change under stabilization.  The second
part, which involves the analysis of a tame isomorphism $\Psi\co  (\alg,
\df) \to (\alg', \df')$, will become important in
\fullref{ssec:invariance}, so a rough sketch of the proof is in
order.  Let an augmentation $\varepsilon'$ on $(\alg', \df')$ be
given.  The goal is to find a corresponding augmentation $\varepsilon$
on $(\alg, \df)$ and an isomorphism between the chain complexes $(A,
\dfe_1)$ and $(A', (\df')^{\varepsilon'}_1)$.  Chekanov proved that
any tame isomorphism, such as $\Phi^{\varepsilon'} \circ \Psi$, can be
factored into a composition $\overline{\Psi} \circ \Psi^0$, where the image
of a generator $q$ under $\overline{\Psi}$ has no constant terms and
$\Psi^0$ is of the form $\Psi^0(q_i) = q_i+c_i$ for some $c \in
\zz/2$.  It is possible to prove that the constants $c_i$ come from an
augmentation $\varepsilon$ of $(\alg, \df)$, and that $\overline{\Psi}$
conjugates $\df^\varepsilon$ and $(\df')^{\varepsilon'}$.  Since
$\overline{\Psi}$ is nondecreasing in the tensor powers of $\alg$, it
restricts to a chain isomorphism between the linearized complexes $(A,
\dfe_1)$ and $(A', (\df')^{\varepsilon'}_1)$.  This argument can
clearly be run in the other direction, thus producing a bijection
between $\mathcal{P}$ and $\mathcal{P}'$.  For more details, see
\cite{chv} or \cite{lecnotes}.

\begin{rem} 
  The set $\mathcal{P}$ of Poincar\'e--Chekanov polynomials is not
  necessarily a one-element set: Melvin and Shrestha
  \cite{melvin-shrestha} found examples of Legendrian knots with
  arbitrarily large sets.
\end{rem}

\begin{exa} \label{ex:5-2-lin}
  Referring back to \fullref{ex:5-2-diff}, it is not hard to check
  that there is a unique augmentation of $(\alg, \df)$ in which both
  of the degree $0$ generators $q_1$ and $q_2$ are sent to $1$ by the
  augmentation.  The resulting linearized differential is:
  \begin{equation} \label{eqn:5-2-dfe}
    \begin{split}
      \df^\varepsilon_1 q_1 &= \df^\varepsilon_1 q_2 = 
      \df^\varepsilon_1 q_6 = \df^\varepsilon_1 q_7 = 0 \\
      \df^\varepsilon_1 q_3 &= q_1 + q_2 \\
      \df^\varepsilon_1 q_4 &= q_1 + q_2 \\
      \df^\varepsilon_1 q_5 &= q_3 + q_4 \\
      \df^\varepsilon_1 q_8 &= q_1 \\
      \df^\varepsilon_1 q_9 &= q_1.
    \end{split}
  \end{equation}
  An easy computation shows that the linearized homology is generated
  by $[q_6]$, $[q_7]$, and $[q_8+q_9]$, so the set of
  Poincar\'e--Chekanov polynomials is
  \begin{equation*}
    \bigl\{t^{-2} + t + t^2 \bigr\}.
  \end{equation*}
\end{exa}
Notice that the Poincar\'e polynomial in the example is symmetric
about degree $0$, with the exception of a class in $H_1(A,
\df_1^\varepsilon)$.  \fullref{thm:pd} asserts that this symmetry
holds in general.

\section{The expanded algebra}
\label{sec:expanded-alg}

As suggested in the introduction, the proof of duality for the
linearized DGA requires disks with two positive corners rather than
the single positive corner of the disks in
\fullref{defn:ce-immersed}. One method to access these disks
comes from the following construction: given a Legendrian knot $K$,
shrink it so that its Reeb chords have length at most $\delta>0$.  Let
$f\co  \rr^3 \to \rr^3$ be the vertical translation by $1$, and let
$K^{n}$ be the link with components $\{K, f(K), \ldots, f^{n-1}(K)\}$.
The Reeb chords for this link are far from isolated, but after a
perturbation using a Morse function on the knot $K$, the chords fall
into three families:
\renewcommand{\descriptionlabel}[1]{\hspace{\labelsep}\textbf{#1}}
\begin{description}
\item[$q$ chords\ ] Chords that start on translates of the bottom
  strand of a chord for $K$ and end on translates of the top strand
  (this includes the original chords of $K$),
\item[$p$ chords\ ] Chords that start on translates of the top strand
  of a chord for $K$ and end on translates of the bottom strand, and
\item[$c$ and $d$ chords\ ] Chords corresponding to critical points of
  the Morse function.
\end{description}
The perturbation procedure and the families of chords will be
described in more detail in Sections~\ref{ssec:expanded-gens} and
\ref{ssec:expanded-diffl}.  After identifying chords that differ by a
translation, these chords will generate the $n$--fold expanded algebra
$\falg^{n}$ of the knot.  For the proof of duality, only $\falg^2$
will be necessary. The higher-order algebras can be used to define
higher cohomology operations; for example, $\falg^3$ is required to
define the inverse of the duality map.

The $p$ chords have the following useful property: a disk with a
negative corner at a $p$ chord locally projects to the same quadrant
as a disk that has a positive corner at a $q$ chord. In other words,
an output of $p$ is the same as an input of $q$. This indicates that,
at least at the linear level, $p$ chords could be used to represent
generators for the cochain complex of the original knot.  It follows
that a disk in the Lagrangian diagram of $K^{n}$ with, say, a positive
corner at a $q$ chord and a negative corner at a $p$ chord corresponds
to a disk with two positive corners in the Lagrangian diagram of the
original knot $K$; see \fullref{fig:reinterpret-disk-1}.  This
property will be exploited to define the duality map in
\fullref{sec:duality}.

\begin{rem}
  This construction was extracted from embedding a ``small''
  Legendrian knot in the contact circle bundle $\rr^2 \times S^1$,
  which is equivalent to taking infinitely many vertical translates of
  the original knot in $\rr^3$.  The full ``Morse--Bott'' theory for
  Legendrian knots in circle bundles, which inspired the language and
  notation of this section, was fully worked out by Sabloff in
  \cite{s1bundles}. See also the last section of Ekholm, Etnyre and
  Sullivan \cite{ees:ori} for a similar construction.
\end{rem}

\subsection{Generators of the algebra}
\label{ssec:expanded-gens}

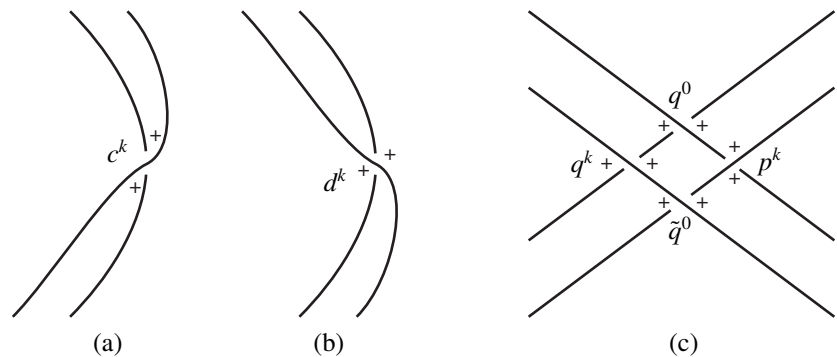
\begin{figure}[ht!]\small
\labellist
\pinlabel $c^k$ <2pt, 0pt> at 38 63
\pinlabel $d^k$ <-2pt, 0pt> at 124 52
\pinlabel $p^k$ <2pt, 2pt> at 284 57
\pinlabel $q^k$ at 215 58
\pinlabel $q^0$ at 252 84
\pinlabel {$\wtilde q^0$} <0pt, -2pt> at 252 35
\endlabellist
  \centerline{\includegraphics{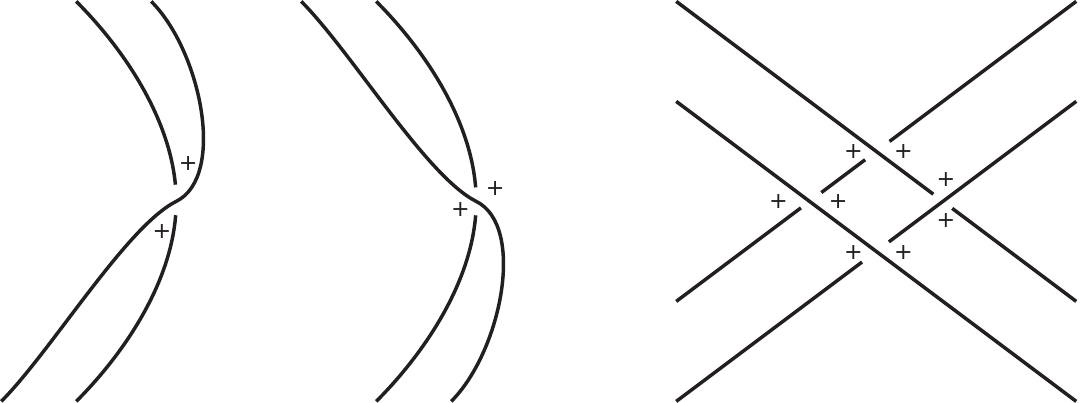}}
  \cl{(a)\hglue 1in(b)\hglue1.7in(c)\hglue.3in}  
\caption{The effect of the perturbation on the crossings of $K$ and
    $f^k(K)$ near (a) a maximum, (b) a minimum, and (c) a crossing of
    $K$.  For this particular choice of $g$, the critical points of
    $\tilde{g}|_K$ to the left of the crossing are maxima and are
    minima to the right.  A different choice of $g$ would give a
    different configuration of signs and of the position of the
    surrounding maxima and minima.}
  \label{fig:perturb}
\end{figure}

In order to define the expanded algebra based on the
Chekanov--Eliashberg DGA of $K^n$, a more precise description of the
perturbation is necessary.  Let $g$ be a function on the $xy$ plane,
supported in a neighborhood of the Lagrangian diagram $\pi_l(K)$, and
let $\tilde{g}$ be a $z$--invariant lift to $\rr^3$.  Choose $g$ so
that it has no critical points in a neighborhood of $\pi_l(K)$ and so
that $\tilde{g}|_K$ is a Morse function whose critical points do not
lie above the crossings of $\pi_l(K)$.  The perturbation of $K^n$ is a
sequence of small, but progressively larger, shifts of $f^j(K)$ in the
direction $i \nabla g$, where $i$ is the usual complex structure on
$\rr^2$.  \fullref{fig:perturb} shows that effects of this
perturbation near the critical points and crossings.

As indicated above, the perturbation introduces a host of new
generators.  For example, \fullref{fig:perturb} shows that every
critical point of $\tilde{g}|_K$ gives rise to a new chord between any
two components of $K^n$.  Label the maxima of $\tilde{g}|_K$ by
$\{c_1, \ldots, c_r\}$ and the minima by $\{d_1, \ldots, d_r\}$.  At a
crossing of $K$, there are four chords for every pair of components:
the original chord ($q^0$ in the figure), its vertical translate
($\tilde{q}^0$), a chord from the bottom strand of $K$ to the top
strand of $f^k(K)$ ($q^k$), and a chord from the top strand of $K$ to
the bottom strand of $f^k(K)$ ($p^k$). In fact, each critical point
and crossing of $K$ gives rise to a half-lattice or lattice of
crossings of $K^n$.  In the expanded algebra, unlike the
Chekanov--Eliashberg DGA of $K^n$, the original chord $q^0$ and its
translate $\tilde{q}^0$ will be identified.  In general, \emph{the
  expanded algebra identifies any two chords that, modulo
  perturbation, are vertical translates of one another.}  This is the
only difference between the expanded algebra and the full
Chekanov--Eliashberg DGA of $K^n$.

More formally, the algebra is defined as follows:

\begin{defn}
  The \textit{$n$--fold expanded algebra} $\falg^n$ is the graded free
  unital tensor algebra over $\zz/2$ generated by the following:
  \begin{itemize}
  \item To the $i^{th}$ double point, associate generators $\{q_i^k,
    p_i^k\}_{k=0,1, \ldots, n-1}$. The generator $q^k_i$
    (respectively, $p^k_i$) corresponds to the chord that starts on
    the bottom (resp.\ top) strand of $K$ and ends on the top
    (resp.\ bottom) strand of $f^k(K)$.\footnote{The generators $p^k$
      would be labeled as $p^{k-1}$ in \cite{s1bundles}; the slight
      change in notation greatly improves readability in the current
      context.} The number $k$ is the \textit{level} of the generator.
  
  \item To each point $c_i$, associate generators $\{c_i^k\}_{k=1, 2,
      \ldots, n-1}$, and to each point $d_i$, associate generators
    $\{d_i^k\}_{k=1, 2, \ldots, n-1}$.  The generators $c^k_i$ (resp.\
    $d^k_i$) represent the Reeb chords that start at a maximum
    (resp.\ minimum) of $\tilde{g}|_K$ and end
    at the corresponding point on $f^k(K)$.
  \end{itemize}
\end{defn}

Note that, by construction, we have:
\begin{equation} \label{eqn:heights}
  \begin{aligned}
    h(q_i^k) &> k, & h(c_i^k) &\approx k, \\
    h(p_i^k) &< k, & h(d_i^k) &\approx k.
  \end{aligned}
\end{equation}
The gradings of the new generators are easy to define.  Let $q^0_i$
inherit its grading from the DGA $(\alg, \df)$; the others are graded
as follows:
\begin{equation} \label{eqn:gradings}
  \begin{aligned}
    |q_i^k| &= |q_i^0|, & |c_i^k| &= 0, \\
    |p_i^k| &= -1-|q_i^0|, & |d_i^k| &= -1.
  \end{aligned}
\end{equation}

\subsection{The differential}
\label{ssec:expanded-diffl}

The definition of the differential \fdf\ on the $n$--fold expanded
algebra is straightforward: it is simply the contact homology
differential for the link $K^n$, up to identifications of chords that
differ by a vertical translation.  To show that the differential is
well-defined, it is useful to introduce Mishachev's relative homotopy
splitting of a link DGA \cite{kirill}, as interpreted by Ng
\cite{lenny:computable}.  The motivation for this language is the
following simple observation: while traversing counter-clockwise the
boundary of an immersed disk involved in the definition of \df, the
upper strand of the first negative corner and the lower strand of the
last negative corner agree with the upper and lower strands,
respectively, of the positive corner.  On the other hand, the lower
strand of the first of two consecutive negative corners coincides with
the upper strand of the second.

More formally, label the components of $K^n$ from bottom to top by
$\{0, \ldots, n-1\}$.   For $j \neq k$, let
$\Gamma_{jk}$ be the submodule of \falg\ generated by words of the
form $x_1 \cdots x_m$ so that:
\begin{itemize}
\item The upper strand of $x_1$ is $j$,
\item The lower strand of $x_m$ is $k$, and 
\item The lower strand of $x_i$ coincides with the upper strand of
  $x_{i+1}$ for $1 \leq i \leq m-1$.
\end{itemize}
If $j=k$, include an indeterminate $e_j$ as a generator of
$\Gamma_{jj}$ as well.  For $K^n$, the component $j$ lies above the
component $k$ if and only if $j > k$, so it is easy to check that
$\Gamma_{jk}$ is empty for $j<k$ and that for $j \geq k$, the upper
(resp.\ lower) strand of $x_i$ in a generator $x_1 \cdots x_m$ is at
most $j$ (resp.\ at least $k$).

The key fact about this construction, whose proof is essentially
outlined above, is that the module $\Gamma_{jk}$ is closed under the
action of the differential \df:

\begin{lem}{\rm \cite[Lemma 2.19]{lenny:computable}}\label{lem:splitting}\qua
  If the upper strand of a chord $x$ is $j$ and its lower strand is
  $k$, then $\df x \in \Gamma_{jk}$, where a $1$ in the differential
  of $x \in \Gamma_{jj}$ is replaced by $e_j$.
\end{lem}

More is true: the differential module $(\Gamma_{jk}, \df)$ is
invariant under Legendrian isotopy up to stable tame isomorphism.
Further, if the set of augmentations is restricted to ones that vanish
on chords between components, then the set of linearized homologies of
$(\Gamma_{jk}, \df)$ is also invariant.

To prove that \fdf\ is well-defined, suppose $x$ is an arbitrary generator
for the Chekanov--Eliashberg DGA of $K^n$ and let $f^k(x)$ be a
vertical translate of $x$ (up to perturbation).
\fullref{lem:splitting} and the observation before it regarding the
upper and lower strands of the $x_i$ that make up the generators of
$\Gamma_{jk}$ together imply that any disk contributing to $\df x$
only involves the link consisting of components of $K^n$ that lie
between the components of the upper and lower strands of $x$
(inclusive).  Similarly, any disk contributing to $\df f^k(x)$ only
involves the link lying between the the components of the upper and
lower strands of $f^k(x)$ (inclusive).  The Lagrangian diagrams of
these links are the same (up to a small shift), so the differentials
of $x$ and $f^k(x)$ are the same up to identification of
vertically-translated generators.  In particular, the differential
\fdf\ on \falg\ is well-defined.

Before moving on, a few words are in order about the form of the
differential.  Following the language of
\cite{kirill},\footnote{Despite the similarities in language and
  diagrammatics, the $N$--copy of \cite{kirill} and the perturbed link
  $K^n$ are quite different objects.  In particular, $K^n$ is
  unlinked.}  define a smooth ``stick-together'' map $s\co  \rr^2 \to
\rr^2$ that retracts the Lagrangian diagram of the perturbed link
$K^n$ onto the Lagrangian diagram of $K$.  The disks that define the
differential split into two classes: the first class consists of
\textit{thick disks} that the stick-together maps sends to immersed
disks in the original diagram.  As can be seen in
\fullref{fig:perturb}, a disk in $K^n$ with a positive corner at a
$q$ crossing will have a positive corner in the original diagram,
while a disk in $K^n$ with a positive $p$ corner will have a negative
corner in the original diagram; a parallel correspondence holds at
negative corners.

The second class of disks consists of \textit{thin disks} whose images
under the stick-together map lie in the original diagram of $K$.  A
thin disk has one of the following forms, as can be seen from
examining \fullref{fig:perturb}:
\begin{enumerate}
\item A disk that flows down from a positive corner at a $c$ crossing to a
  negative corner at one of the two adjacent $d$ crossings,
\item A disk that flows down from a positive corner at a $c$ crossing
  to negative $p$ and $q$ corners at a crossing of $\pi_l(K)$ that
  lies between the $c$ and an adjacent $d$,
\item A disk that flows down from a positive corner at a $q$
  (resp.\ $p$) crossing and a negative corner at another $q$
  (resp.\ $p$) crossing in the same lattice to a negative corner at the
  next $d$ crossing,\footnote{More than the two copies of $K$ depicted
    in \fullref{fig:perturb} are necessary to see the disk with
    $p$ corners.}
\item A disk that flows down from a positive corner at a $q$
  (resp.\ $p$) crossing and a negative corner at another $q$
  (resp.\ $p$) crossing in the same lattice to negative $p$ and $q$
  corners in another lattice (of course, there can be no intervening
  $c$ or $d$ crossings), or
\item A disk that lies entirely inside the lattice of crossings
  created by the perturbation.
\end{enumerate}
These disks have a close relationship with the negative gradient flow
of $\tilde{g}|_K$: disks of type (1) correspond to flowlines from
maximum to minimum, while disks of types (2--4) coincide with
``partial flowlines.''  In each case, the positive corner lies at the
end of the flowline or partial flowline with the largest value of $g$
and flows down to a negative corner (possibly with other negative
corners at each end).

It is easy to see that the differential of a $c$ generator comes
entirely from thin disks of types (1) and (2), while the differential
of a $d$ generator comes entirely from thick disks and thin disks of
type (5).  There are exactly two thin disks of type (3) that
contribute to the differential of each $q$ and each $p$ generator, and
possibly more of types (4) or (5).

\begin{figure}[ht!]
  \centerline{\includegraphics[scale=.9]{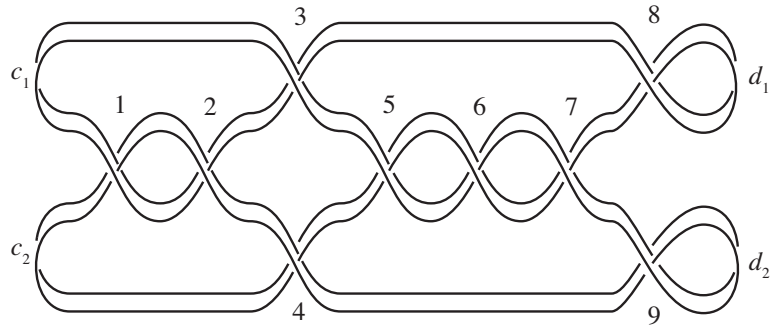}}
  \caption{A perturbation of the knot in \fullref{fig:numbered-5-2}
    that gives a $2$--fold expanded algebra.}
  \label{fig:perturb-5-2}
\end{figure}

\begin{exa} \label{ex:5-2-expanded}
  Consider the $2$--fold expanded algebra of the knot in
  \fullref{fig:numbered-5-2} given by the perturbation in
  \fullref{fig:perturb-5-2} (the perturbation comes from shifting
  the entire diagram down).  Consider the differential of $q^1_6$.
  The thick disks contribute the following terms:
  \begin{equation}
    p_5^1 + p_7^1.
  \end{equation}
  Thin disks of types (3)--(5) contribute quite a few more:
  \begin{equation}
    q_6^0(q_7^0 p_7^1 + q_9^0 p_9^1 + d_2)
    + (p_7^1 q_7^0 + p_8^1 q_8^0 + d_1)q_6^0
    + q_6^0 p_6^1 q_6^0.
  \end{equation}
\end{exa}

\begin{figure}[ht!]
  \centerline{\includegraphics[scale=.8]{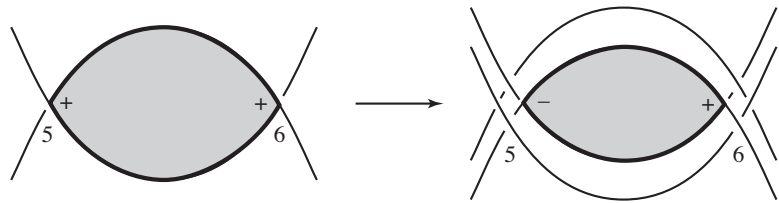}}
  \caption{A disk with two positive corners in the Lagrangian diagram
    of the original knot $K$ corresponds to a disk with a positive
    corner at a $q$ chord and a negative corner at a $p$ chord in a
    Lagrangian diagram of $K^2$.}
  \label{fig:reinterpret-disk-1}
\end{figure}

To connect this example to the upcoming proof of duality, recall the
observation at the beginning of this section that a disk in the
Lagrangian diagram of $K^{n}$ with, say, a positive corner at a $q$
chord and a negative corner at a $p$ chord corresponds to a disk with
two positive corners in the Lagrangian diagram of the original knot
$K$.  \mbox{\fullref{fig:reinterpret-disk-1}} shows this explicitly in the
example above: a disk with two positive corners in the Lagrangian
diagram of the original knot $K$ yields a $p^1_5$ term in $\fdf q^1_6$
in the expanded algebra of $K^2$.

\section{The proof of duality}
\label{sec:duality}

When linearized, the expanded algebra carries a lot of structure.  In
this section, this structure will be uncovered and exploited to prove
\fullref{thm:pd}.

\subsection{Structure of the linearized expanded algebra}
\label{ssec:linear-structure}

This section examines the following decomposition of the linear pieces
of $\falg^n$: let $Q^m$ be the vector space generated by the level $m$
chords $\{ q_i^m, c_i^m, d_i^m\}$ and let $P^m$ be the vector space
generated by $\{p_i^m\}$.  The goal of this section is to understand
the relationships between the $Q$ and $P$ vector spaces once $\falg$
has been linearized; the primary tools will be Lemmas~\ref{lem:stokes}
and \ref{lem:splitting}.

\subsubsection{Extending the Augmentation}

The first step will be to extend an augmentation for $\alg$ to all of
$\falg^n$.  Given an augmentation $\varepsilon$ of \alg, there is a
trivial extension to a map $\widehat{\varepsilon}$ on $\falg^n$:
simply set $\widehat{\varepsilon}$ to be zero on all generators not in
$Q^0$.  This map turns out to be an augmentation for the expanded
algebra:

\begin{prop} \label{prop:ext-augm}
  $\widehat{\varepsilon} \circ \fdf = 0$.
\end{prop}

\begin{proof}
  The first step is to show that every term in the differential of a
  generator not in $Q^0$ contains at least one generator not in $Q^0$.
  This is simple: since a generator $x^m$, $m \geq 1$, corresponds to
  a crossing between two different components of $K^n$,
  \fullref{lem:splitting} implies that there must be at least one
  generator in each term of $\fdf x^m$ whose upper and lower strands
  lie on different components of $K^n$.  This implies that
  $\widehat{\varepsilon} \circ \fdf x^m = 0$ for $m \geq 1$ since
  $\widehat{\varepsilon}$ vanishes on generators of level $m \geq 1$.

  This leaves $Q^0$.  When restricted to $Q^0$, it turns out that the
  expanded differential \fdf\ agrees with the original differential
  \df\ on $K$.  Let $q \in Q^0$. Since the Reeb chords that generate
  $Q^0$ are shorter than any other chords, \fullref{lem:stokes}
  implies that $\fdf q$ contains only other generators from $Q^0$. In
  particular, the boundary of any disk contributing to $\fdf q$ lies
  on a single component of $K^n$, and hence contributes to the
  original differential of $q$ on $K$. Thus, on $Q^0$, $\fdf = \df$
  and the proposition follows from the definition of
  $\widehat{\varepsilon}$.
\end{proof}

\subsubsection{Splitting the linearized differential}

The linearized differential $\fdf^{\widehat{\varepsilon}}_1$ has a
rich structure that will be exploited in the proof of duality.  For
convenience, the superscript $\widehat{\varepsilon}$ will be dropped
from the notation henceforth, ie, all differentials will be assumed
to be augmented.  The first important property of the linearized
differential is that it preserves level:

\begin{lem} \label{lem:preserve-level} The images of the linearized
  differential obey $\fdf_1(Q^0) \subset Q^0$ and\break $\fdf_1(Q^m \oplus
  P^m) \subset Q^m \oplus P^m$ for all $m \geq 1$.  Further,
  $\fdf_1(P^m) \subset P^m$.
\end{lem}

\begin{proof}
  A summand $y_1 \cdots y_k$ of $\fdf x^m$ will contribute $y_i$ to
  the linearized differential if $\widehat{\varepsilon}(y_j) = 1$ for
  all $j \neq i$.  By the definition of $\widehat{\varepsilon}$, all
  such $y_j$ must be in $Q^0$.  Hence, the only place that the
  boundary of the disk giving rise to the summand $y_1 \cdots y_k$ can
  switch from the upper strand of $x^m$ to the lower strand is at
  $y_i$.  Thus, $y_i$ must have the same upper and lower strands as
  $x^m$, so the image of the linearized differential of $x^m$ must
  also have level $m$.
  
  The second statement follows from the first, \fullref{lem:stokes},
  and Equation \eqref{eqn:heights}, which show that a term in $\fdf_1
  p^m$ must have level $m$, have height less than $p^m$, and hence
  must be an element of $P^m$.
\end{proof}

On any given level $m \geq 1$, the linearized differential splits into
three pieces by restricting domains and ranges:
\begin{enumerate}
\item $\fdf_Q\co  Q^m \to Q^m$,
\item $\fdf_P\co  P^m \to P^m$, and
\item $\eta\co  Q^m \to P^m$.  This is the map that appears schematically
  in \fullref{fig:reinterpret-disk-1}.
\end{enumerate}

That there is no piece of the linearized differential with domain
$P^m$ and range $Q^m$ is a consequence of the previous lemma. 

\begin{exa}
  As noted in \fullref{ex:5-2-lin}, the only generators augmented
  in the case of the knot in \fullref{fig:numbered-5-2} are $q_1^0$
  and $q_2^0$.  Thus, for the generators $q^1_6$ in
  \fullref{fig:perturb-5-2}, the thin disks computed in
  \fullref{ex:5-2-expanded} contribute nothing to the linearized
  differential.  Thus, it is straightforward to read off that $\fdf_Q
  (q_6^1) = 0$.  The $\eta$ map at this crossing, however, is
  nontrivial:
  \begin{equation} \label{eqn:5-2-eta} \eta \left( q_6^1 \right) =
    p_5^{1} + p_7^{1}.
  \end{equation}
\end{exa}

The components of the linearized differential give rise to two
families of chain complexes:

\begin{prop} \label{prop:decomposition}
  \begin{enumerate}
  \item $(Q^m,\fdf_Q)$ and $(P^m,\fdf_P)$ are chain complexes, and
  \item $\eta\co  (Q^m, \fdf_Q) \to (P^m, \fdf_P)$ is a chain map of
    degree $-1$ for $m \geq 1$.
  \end{enumerate}
\end{prop}

\begin{proof}
  The fact that the degrees of $\fdf_Q$, $\fdf_P$, and $\eta$ are all
  $-1$ follows from the fact that they are components of the degree
  $-1$ differential.  

  To prove the remainder of the proposition, represent $\fdf_1\co  Q^m
  \oplus P^m \to Q^m \oplus P^m$ by the following matrix:
  \begin{equation*}
    \begin{bmatrix}
      \fdf_Q & 0 \\
      \eta & \fdf_P
    \end{bmatrix}.
  \end{equation*}
  Since $\fdf^2 = 0$, every component of the following matrix is also
  zero:
  \begin{equation*}
    \begin{bmatrix}   \fdf_Q & 0 \\  \eta & \fdf_P \end{bmatrix}^2 =
    \begin{bmatrix} \fdf_Q^2 & 0 \\ \eta \fdf_Q + \fdf_P \eta &
      \fdf_P^2 \end{bmatrix}.
  \end{equation*}
  This finishes the proof.
\end{proof}

The proof of \fullref{prop:ext-augm} shows that one of the new
chain complexes in the proposition above is familiar:

\begin{cor} \label{cor:a-q0}
  $(A, \df_1) \simeq (Q^0, \fdf_Q)$
\end{cor}

Here is one more corollary of \fullref{prop:decomposition}:

\begin{cor} \label{cor:mapping-cone} For $m \geq 1$, $(Q^m \oplus P^m,
  \fdf_1)$ is the mapping cone for the chain map $\eta$.
\end{cor}

The differential $\fdf_Q$ splits further in a manner similar to the
original $\fdf_1$.  For $m \geq 1$, let $\overline{Q}^m$ be the vector
space generated by $\{ q_i^m \}$ and let $C^m$ be generated by
$\{c_i^m, d_i^m\}$. As above, split $\fdf_Q$ by restricting domains
and ranges:
\begin{enumerate}
\item $\fdf_{\overline{Q}}\co  \overline{Q}^m \to \overline{Q}^m$,
\item $\fdf_C\co  C^m \to C^m$, and
\item $\rho\co  \overline{Q}^m \to C^m$.
\end{enumerate}
By the same argument as in the proof of
\fullref{lem:preserve-level}, there is no component of $\fdf_Q$ with
domain $C^m$ and range $Q^m$.  As in
\fullref{prop:decomposition}, it follows that $(\overline{Q}^m,
\fdf_{\overline{Q}})$ and $(C^m, \fdf_C)$ are both chain complexes and
that $\rho$ is a chain map between them.  Interestingly, though perhaps not
surprisingly given the Morse--Bott motivation for the circle bundle
theory from which the expanded algebra is derived, both of these chain
complexes are familiar:

\begin{lem} \label{lem:q-bar}
  $(\overline{Q}^m, \fdf_{\overline{Q}}) \simeq (Q^0, \fdf_Q)$
\end{lem}

\begin{proof}
  There is an obvious correspondence between generators and their
  degrees, so it suffices to extend this correspondence to the
  differentials.  As demonstrated in the proof of
  \fullref{lem:preserve-level}, a disk contributing to
  $\fdf_{\overline{Q}} q^m$ must have a positive corner at $q^m$ and
  negative corners at a $q$ crossing of level $m$ and possibly others
  at augmented $q^0$ corners. Thus, under the stick-together map, a
  thick disk that contributes to $\fdf_{\overline{Q}}q^m$ becomes a
  disk contributing to $\fdf_Q q^0$.  Further, no thin disk can
  contribute to $\fdf_{\overline{Q}}$ since any thin disk with a
  positive $q$ corner must have either a negative corner at a $d$
  (type (3)) or at a $p$ (types (4) and (5)).  Overall, then, the
  differentials on $\overline{Q}^m$ and on $Q^0$ are constructed out
  of the same disks, and hence are the same under the correspondence
  between generators.
\end{proof}

\begin{lem} \label{lem:morse-witten} Up to shifting degree down by
  $1$, $(C^m, \fdf_C)$ is isomorphic to the Morse--Witten complex
  $(C^{MW}, \df^{MW})$ of the circle with respect to $\tilde{g}|_K$.
\end{lem}

\begin{proof}
  It is clear from the construction that there is a one-to-one
  correspondence between the generators of $C^m$ and those of
  $C^{MW}$.  As noted in \fullref{ssec:expanded-diffl}, the
  differential of a $c$ chord comes completely from thin disks of
  types (1) or (2).  Disks of type (2) have a negative $p$ corner, so
  cannot contribute to $\fdf_C$.  Disks of type (1) are in one-to-one
  correspondence with negative gradient flowlines for $g$, and hence
  $\fdf_C c^m = \df^{MW} c^m$.  For degree reasons, $\fdf_C d^m = 0 =
  \df^{MW} d^m$.  This completes the proof.
\end{proof}

\subsubsection{``Vertical'' relationships}

There is a ``vertical'' relationship between any pair of complexes
$(Q^m, \fdf_Q)$ and $(Q^{m+1}, \fdf_Q)$.  This relationship can be
encoded by a \textit{translation map} $\tau\co  Q^m \to Q^{m+1}$ that
raises the level of a generator by $1$; a similar map exists for
$P^m$.  The maps and spaces defined thus far are pictured in
\fullref{fig:maps}.  The translation map interacts nicely with the
linearized differential:

\begin{figure}[ht!]
  \begin{equation*}
    \xymatrix{
      \vdots & \vdots \\
      Q^2 \ar[r]^\eta \ar[u]^\tau & P^2 \ar[u]_\tau \\
      Q^1 \ar[u]^\tau \ar[r]^\eta & P^1 \ar[u]_\tau  \\
      Q^0 \ar[u]^\tau & 
    }
  \end{equation*}
  \caption{Relationships between the maps and spaces defined thus far.
    By \fullref{lem:translation}, the diagram commutes.}
  \label{fig:maps}
\end{figure}
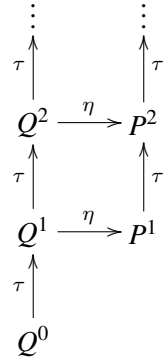

\begin{lem}
  \label{lem:translation} The translation map $\tau$ is a chain
  isomorphism on $(P^m, \fdf_P)$ and on $(Q^m, \fdf_Q)$ for all $1
  \leq m \leq n-2$.  Further, $\tau$ commutes with $\eta$.
\end{lem}

\begin{proof}
  For all levels in $P$ and for levels greater than zero in $Q$, the
  translation map $\tau$ is obviously bijective.

  To prove that $\tau$ is a chain map and commutes with $\eta$,
  consider \fullref{lem:preserve-level} and the definition of
  $\widehat{\varepsilon}$.  These combine to show that any disk that
  contributes to the linearized differential of a generator of level
  $m$ has one negative corner at a crossing between $K$ and $f^m(K)$
  and all other negative corners at (augmented) self-crossings of $K$
  or of $f^m(K)$.  In particular, such disks depend only on the link
  $K \cup f^m(K)$.  For $m \geq 1$, the Lagrangian projections of the
  links $K \cup f^m(K)$ and $K \cup f^{m+1}(K)$ are combinatorially
  identical, and hence so are the disks that contribute to $\fdf_1$
  --- and hence to $\fdf_Q$, $\fdf_P$, and $\eta$ --- at levels $m$ and
  $m+1$.
\end{proof}

More interesting phenomena occur when comparing $(Q^0, \fdf_Q)$ to
$(Q^1, \fdf_Q)$. The chain complex $(Q^1, \fdf_Q)$ can be written as
follows, with subscripts denoting degrees:
\begin{equation*}
  \xymatrix{
    \ar[r]^{\fdf_{\overline{Q}}} & \overline{Q}^1_2 \ar[r]^{\fdf_{\overline{Q}}} & 
    \overline{Q}^1_1 \ar[r]^{\fdf_{\overline{Q}}} \ar[dr]^\rho & 
    \overline{Q}^1_0 \ar[r]^{\fdf_{\overline{Q}}} \ar[dr]^\rho & 
    \overline{Q}^1_{-1} \ar[r]^{\fdf_{\overline{Q}}} & 
    \overline{Q}^1_{-2} \ar[r]^{\fdf_{\overline{Q}}} & \\
    & & 0 \ar[r] & C^1_0 \ar[r]^{\fdf_C} & C^1_{-1} \ar[r] & 0 &
  }
\end{equation*}
It follows from \fullref{lem:q-bar} that:
\begin{equation} \label{eqn:h-k} \dim H_{k}(Q^1, \fdf_Q) = \dim
  H_k(\overline{Q}^1, \fdf_{\overline{Q}}) = \dim H_{k}(Q^0, \fdf_Q)
\end{equation}
for $k>1$ and $k<-1$.  This, together with knowledge of the case
$k=-1$, will suffice for the upcoming proof of duality.

In degree $-1$, it turns out that
\begin{equation} \label{eqn:h-1}
  \dim H_{-1}(Q^1, \fdf_Q) = \dim H_{-1}(Q^0, \fdf_Q) + 1.
\end{equation}
On one hand, the kernel of $\fdf_Q$ on $Q^1_{-1}$ is the direct sum of
the kernel of $\fdf_{\overline{Q}}$ on $\overline{Q}^1_{-1}$ and all
of $C^1_{-1}$, the space generated by the $d$'s. On the other, the
image of $\fdf_Q$ on $Q^1_0$ comes from studying the images of
$\fdf_{\overline{Q}}$, $\fdf_C$, and $\rho$.
\fullref{lem:morse-witten} implies that the $c^1$ terms have a
linearized differential given by
\begin{equation} \label{eqn:c}
  \fdf_Q (c^1) = d^1 + \overline{d}^1,
\end{equation}
where $d$ and $\overline{d}$ are adjacent to $c$ along the knot.
Because of the configuration of positive and negative corners around a
$d$ crossing (see \fullref{fig:perturb}), only thin disks of type
(3) that begin at an augmented crossing contribute to the image of
$\rho$.  In fact, by the discussion of how thin disks contribute to
the differential $\fdf$ in \fullref{ssec:expanded-diffl}, if $q^0_i$ is
augmented, then $\rho(q^1_i)$ is the sum of the two $d^1$ generators
adjacent to $q_i$.  This is also the image under $\fdf_C$ of the sum
of the $c^1$ generators that lie between the two $d^1$ generators on
the knot $K$.\footnote{There are two such sets of $c^1$ generators,
  depending on the direction in which the knot $K$ is
  traversed. Either set works.  This is true even if one is the empty
  set, as this is the case that the two adjacent $d^1$ generators
  coincide and hence sum to $0$ over $\zz/2$.} Thus, it makes sense to
consider the following change of basis for $Q^1_0$: let $\beta_i$ be
the sum of the $c^1$ generators that lie between the two $d^1$
generators whose sum is $\rho (q^1_i)$.  Change basis so that if
$q^0_i$ is augmented, then $q^1_i$ is replaced by $q^1_i + \beta_i$.
In this basis, it is clear that the image of $\fdf_Q$ is the direct
sum of the image of $\fdf_{\overline{Q}}$ on $\overline{Q}^1_0$ and
the image of $\fdf_C$.  Thus, the degree $-1$ homology of $Q^1$ is the
direct sum $H_{-1}(\overline{Q}^1, \fdf_{\overline{Q}}) \oplus
H_0(S^1)$, with the latter summand generated by $[d^1_1] = [d^1_2] =
\cdots = [d^1_r]$.  Equation \eqref{eqn:h-1} follows.

\subsubsection{``Horizontal'' relationships}

There is also a ``horizontal'' relationship between the chain
complexes $(Q^0,\fdf_Q)$ and $(P^1,\fdf_P)$.  Define a degree $1$
pairing\footnote{This assumes that the base ring lies in degree $0$.}
on $P^1 \otimes Q^0$ by:
\begin{equation*}
  \langle p_i,q_j \rangle = \delta_{ij}.
\end{equation*}
This pairing allows us to view $(Q^0,\fdf_Q)$ and $(P^1,\fdf_P)$ as
dual complexes since a disk that contributes $q_j$ to $\fdf_Q (q_i)$
also contributes $p_i$ to $\fdf_P (p_j)$.  Thus, the pairing descends
to homology. In particular, by \eqref{eqn:gradings}:

\begin{lem} \label{lem:hom-dual}
  $\dim H_k(Q^0, \fdf_Q) = \dim H_{-k-1}(P^1, \fdf_P).$
\end{lem}

\begin{exa}
  Combining the lemma with equation \eqref{eqn:5-2-dfe} shows that
  $p_7^1$ represents a non-trivial homology class and that $p_5^1$ is
  a boundary in $(P^1, \fdf_P)$.  Thus, since $\eta$ is a chain map,
  \eqref{eqn:5-2-eta} becomes:
  \begin{equation} \label{eqn:5-2-hom-eta}
    \eta_* [q_6^1] = [p_7^{1}].
  \end{equation}
  Similarly, it is possible to compute that:
  \begin{equation} \label{eqn:5-2-hom-eta-2}
    \eta_* [q_7^1] = [p_6^1].
  \end{equation}
  Since $q_6$ has degree $-2$ and $q_7$ has degree $2$, these
  equations indicate that the $\eta$ map could be the duality map that
  explains the symmetry of the Poincar\'e--Chekanov polynomial for the
  knot in \fullref{fig:numbered-5-2}.
\end{exa}

\subsection{Proof of duality}
\label{ssec:proof}

The analogy between contact homology and Morse theory described in
\fullref{ssec:morse} suggests that the duality map for linearized
Legendrian contact homology should come from disks with two positive
corners. The structure analyzed in the previous section provides the
setting in which to make this analogy precise: disks with two positive
corners are the defining disks for the $\eta$ map, as shown in
\fullref{fig:reinterpret-disk-1}.  The following proposition
demonstrates that the analogy holds in this case:

\begin{prop} \label{prop:eta-isom} The map $\eta\co  (Q^m, \fdf_Q) \to
  (P^m, \fdf_P)$ is an isomorphism.
\end{prop}

\begin{proof}
  By \fullref{cor:mapping-cone}, it suffices to show that $(Q^m
  \oplus P^m, \fdf_1)$ is acyclic.  This, in turn, is obvious, since
  this complex is the same as $\Gamma_{0m}$ for the unlink
  consisting of $K$ and the large vertical translate $f^m(K)$.
\end{proof}

The proof of duality is now a matter of combining
\fullref{prop:eta-isom} with the relationships between the
homologies of $(Q^1, \fdf_Q)$ and $(Q^0, \fdf_Q)$ that were calculated
in \fullref{ssec:linear-structure}.

\begin{proof}[Proof of \fullref{thm:pd}]
  First suppose that $|k|>1$:
  \begin{align*}
    \dim H_{k}(A, \df_1) &= \dim H_{-k-1}(P^1, \fdf_P) && \text{by
      \fullref{lem:hom-dual}} \\
    &= \dim H_{-k}(Q^1, \fdf_Q) && \text{by
      \fullref{prop:eta-isom}} \\
    &= \dim H_{-k}(A, \df_1) && \text{by \eqref{eqn:h-k}}
  \end{align*}
  The proof of the case where $k=1$ is exactly the same up until the
  last step, where \eqref{eqn:h-1} implies that
  \begin{equation*}
    \dim H_{1}(A, \df_1) = \dim H_{-1}(Q^1, \fdf_Q) = \dim H_{-1}(A,
    \df_1) +1, 
  \end{equation*}
  as required by the theorem.  Finally, $\dim H_{0} (A, \df_1)$ does
  not matter for \fullref{thm:pd}, so the proof is complete.
\end{proof}

\subsection{Invariance}
\label{ssec:invariance}

The duality map $\eta$ was defined using a fixed Lagrangian diagram
and perturbing Morse function.  Two pairs of Lagrangian diagrams and
perturbing functions are related by Reidemeister moves, and it is
indeed the case that the duality maps before and after the induced
stable tame isomorphism are, in some suitable sense, conjugate. To set
notation for a precise statement of invariance, let $(\falg, \fdf)$
and $(\falg', \fdf')$ be the expanded algebras of the perturbed
Lagrangian diagrams before and after the Reidemeister moves.  Assume
that both algebras have been appropriately stabilized and that $\Psi$
is a tame isomorphism between them.  Let $\widehat{\varepsilon}$ and
$\widehat{\varepsilon'}$ be augmentations that correspond under $\Psi$
in the manner outlined in \fullref{ssec:linear}; recall that there
is a map $\overline{\Psi}$ that restricts to a chain isomorphism
between the two \emph{linearized} expanded algebras. Further, by
Mishachev's work, $\overline{\Psi}$ must respect levels. Finally,
write the restriction of $\overline{\Psi}$ to $Q^1 \oplus P^1$ in
terms of its components as follows:
\begin{equation}
  \overline{\Psi} = \begin{bmatrix}  
    \Psi_Q & G \\
    H & \Psi_P
  \end{bmatrix}.
\end{equation}
With this notation set, the statement of invariance is:

\begin{prop} \label{prop:invariance} The maps $\Psi_Q$ and $\Psi_P$
  are chain isomorphisms that, on homology, conjugate $\eta_*$ and
  $\eta_*'$.
\end{prop}

\begin{proof}
  The key to the proof is the fact that the map $G\co  P^1 \to (Q')^1$ is
  zero.  This fact comes from examining the tame isomorphisms used to
  prove invariance under the Reidemeister moves; see \cite{chv,ens}.
  In all cases, the isomorphisms are nondecreasing in height, and
  hence a generator of $(Q')^m$ cannot appear in the image of a
  generator of $P^m$.

  With this fact in hand, it follows that the matrix for
  $\overline{\Psi}$ is lower triangular, which immediately shows that
  $\Psi_Q$ and $\Psi_P$ are isomorphisms.  Further, the fact that
  $\overline{\Psi}$ is a chain map implies the following three
  equations:
  \begin{equation}
    \begin{split}
      \fdf'_Q \Psi_Q &= \Psi_Q \fdf_Q, \\
      \fdf'_P \Psi_P &= \Psi_P \fdf_P, \text{ and} \\
      \eta' \Psi_Q + \Psi_P \eta &= \fdf'_P H + H \fdf_Q.
    \end{split}    
  \end{equation}
  The first equation shows that $\Psi_Q$ is a chain isomorphism
  from $(Q^1, \fdf_Q)$ to $(Q'^1, \fdf')$; a similar statement holds
  for $\Psi_P$.  The last equation shows that, on homology, $\eta'_*
  (\Psi_Q)_* = (\Psi_P)_* \eta_*$, as desired.
\end{proof}

\section{Duality and the cap product}
\label{sec:alg-top}

The expanded algebra defined in the previous sections not only
provides the structure with which to prove duality, but also to begin
to examine the ``algebraic topology'' of the linearized Legendrian
contact homology.  In this section, the proof of duality will be
reinterpreted in terms of the classical ``capping with the fundamental
class'' construction.

\subsection{The fundamental class}
\label{ssec:fc}

In Morse theory, the fundamental class of a compact manifold is
represented by the sum of the maxima of a Morse function.
Analogously, the longest Reeb chords --- for example, those coming
from the rightmost cusps in a resolved front diagram --- should give
the fundamental class in degree $1$ in the linearized DGA.  This is
not quite accurate, however, as the change of variables given by an
augmentation can introduce new terms.  One way to overcome this
difficulty --- and to produce a nontrivial class in degree $1$ --- was
discovered by Chekanov \cite[Section 12]{chv}: closing a long Legendrian
knot induces a one-to-one correspondence between Legendrian isotopy
classes of long Legendrian knots and those of Legendrian knots.
Chekanov also proved that the mapping from the set of
Poincar\'e--Chekanov polynomials for a long Legendrian knot to the set
for its closure defined by
\begin{equation*}
  P(t) \longmapsto P(t) + t
\end{equation*}
is bijective.  In particular, there is always a nontrivial homology
class in degree $1$ for the closure. Chekanov's construction, however,
is not precise enough for the needs of this paper.

The necessary precision is derived from the calculations that
established the relationship between $H_{-1}(Q^0, \fdf_Q)$ and
$H_{-1}(Q^1, \fdf_Q)$ in \fullref{ssec:linear-structure}.  Recall
that the proof of equation \eqref{eqn:h-1} decomposed
$H_{-1}(Q^1,\fdf_Q)$ into $H_{-1}(\overline{Q}^1,\fdf_{\overline{Q}}) \oplus
H_0(S^1)$, with the second factor generated by the class $[d]$.  In
essence, the following theorem asserts that the fundamental class is
dual (with respect to the pairing $\langle, \rangle$) to the image of
the generator of $H_0(S^1)$ under the Poincar\'e duality map $\eta$.

\begin{thm} \label{thm:fc-unique} There exists a unique class $\fc$,
  called the {\upshape fundamental class}, in $H_1(A, \df_1)$ that pairs
  to $1$ with $\eta_*([d])$ and to $0$ on the image of
  $H_{-1}(\overline{Q}^1,\fdf_{\overline{Q}})$ under $\eta_*$.
\end{thm}

\begin{proof}
  The theorem follows immediately from \eqref{eqn:h-1} and
  \fullref{prop:eta-isom}.
\end{proof}

The statement of invariance in \fullref{ssec:invariance} only
shows that the fundamental class is defined up to an isomorphism of
the expanded algebra for a given Lagrangian diagram of $K$.  The
fundamental class, however, lies in the original linearized contact
homology of the diagram, so it would be nice if the restriction of the
map $\eta$ to $\overline{Q}^1_{-1}$ were to be independent of the
perturbation function.  To see that this is indeed true, consider a
$q^1$ generator of degree $-1$.  No thin disk can contribute to
$\eta(q^1)$: a thin disk of type $3$ (resp.\ $4$ or $5$) would have
corners at a $d$ crossing (resp.\ a $p$ crossing) and a $q^0$ crossing
of degree $-1$, and hence could not contribute to the linearized
differential.  Since the thick disks depend only on the original
Lagrangian diagram, the image of $q^1$ under $\eta$ does not depend on
the perturbation.  As the coefficient ring is $\zz/2$, the image of
$H_{-1}(\overline{Q}^1,\fdf_{\overline{Q}})$ determines the
fundamental class.  Thus, the restriction of the map $\eta$ to
$\overline{Q}^1_{-1}$ does not depend on the perturbation and can be
read off of the original Lagrangian diagram using thick disks.

The definition of the fundamental class in \fullref{thm:fc-unique}
captures some of the original Morse-theoretic motivation.  If a
Lagrangian diagram is the resolution of a front diagram, then the
longest chords are those coming from the rightmost cusps; see
\fullref{fig:numbered-5-2}.  Fix a chord coming from a right cusp.
Choose a perturbation so that there is a $d$ chord on the loop to the
right of this crossing.  The image of $d$ under $\eta$ is precisely
the right cusp chord, leading to the conclusion that every rightmost
cusp chord is a summand of any representative of the fundamental
class.

\begin{figure}[ht!]
  \centerline{\includegraphics{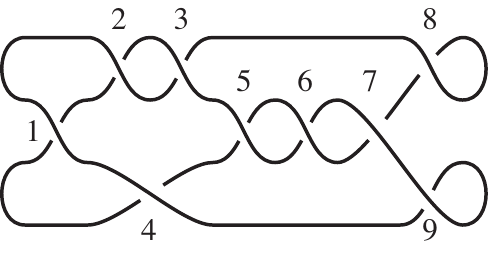}}
  \caption{The resolution of a front diagram of a Figure 8 knot.}
  \label{fig:fig-8}
\end{figure}

\begin{exa}
  The knot in \fullref{fig:numbered-5-2} has a fundamental class
  consisting purely of the right cusps: $[\fc] = [q_8+q_9]$.  This is
  not always the case for resolutions of front diagrams: it is
  possible to calculate that the figure 8 knot in
  \fullref{fig:fig-8} has a fundamental class of $[\fc] =
  [q_5+q_7+q_8+q_9]$ for either of its two augmentations.
\end{exa}

\subsection{The cap product}
\label{ssec:cap}

The goal of this section is to use the cap product and the fundamental
class to exhibit an inverse to the duality map $\eta_*$. The first
step in constructing an inverse to $\eta_*$ is to define the cap
product.  As argued in \fullref{ssec:morse}, the cap product for
linearized Legendrian contact homology should come from disks with one
positive corner and two negative corners.  As with the fundamental
class, the augmentation interferes with this geometric description,
but the length $2$ differential, $\fdf_2$, encodes the relevant
algebraic information.\footnote{Recall that the superscript
  $\hat{\varepsilon}$ has been dropped from the notation for the
  augmented differential.}  Geometrically, the length $2$ differential
counts disks with one positive corner, two negative corners, and
possibly other augmented negative corners.

More technically, the length $2$ differential $\fdf_2$ can be split
into components just as $\fdf_1$ was in
\fullref{ssec:linear-structure}. Write:
\begin{equation}
  P = \bigoplus_{1 \leq m \leq n-2} P^m \quad \text{and} \quad 
  Q = \bigoplus_{0 \leq m \leq n-2} Q^m.
\end{equation}
For now, only the components of $\fdf_2$ with domain $P$ are of
interest; they are:
\begin{enumerate}
\item $\Phi_{QP}\co  P \to Q \otimes P$,
\item $\Phi_{PQ}\co  P \to P \otimes Q$, and
\item $\Phi_{PP}\co  P \to P \otimes P$.
\end{enumerate}
That these preserve total level and are the only components follows
from the same arguments as for \fullref{lem:preserve-level}.  The
first of these maps is a chain map into the tensor product $(Q,
\fdf_Q) \otimes (P, \fdf_P)$.  To see why, look at the length two
component of the equation $\fdf^2 = 0$:
\begin{equation*}
  \fdf_1 \fdf_2 + \fdf_2 \fdf_1 = 0, 
\end{equation*}
and more specifically, the component of the left hand side that maps
$P$ to $Q \otimes P$:
\begin{equation} \label{eqn:df2} \Phi_{QP} \fdf_P + \left( \fdf_Q
    \otimes 1 + 1 \otimes \fdf_P \right) \Phi_{QP} = 0.
\end{equation}
The second map is similarly a chain map.

Pairing the second tensor factor of $\Phi_{QP}$ with $\fc$ turns
this map into the cap product with the fundamental class.  More
precisely:

\begin{defn} \label{defn:eval-cap} Let $k$ be a representative of the
  fundamental class $\fc$. The \textit{evaluation map} $\eval\co  P \to
  \zz/2$ is defined by
  \begin{equation*}
    \eval(p) = \langle p, k \rangle
  \end{equation*}
  on $P^1$ and is extended to the rest of $P$ by $0$.
  
  The \textit{cap product map} $\phi\co  P \to Q$ is defined by
  \begin{equation*}
    \phi = (1 \otimes \eval) \Phi_{QP}.
  \end{equation*}
\end{defn}

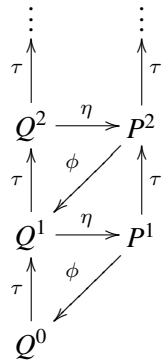
\begin{figure}[ht!]
  \begin{equation*}
    \xymatrix{
      \vdots & \vdots \\
      Q^2 \ar[r]^\eta \ar[u]^\tau & P^2 \ar[u]_\tau \ar[dl]_\phi\\
      Q^1 \ar[u]^\tau \ar[r]^\eta & P^1 \ar[u]_\tau \ar[dl]_\phi \\
      Q^0 \ar[u]^\tau & 
    }
  \end{equation*}
  \caption{Relationships between the maps and spaces defined in
    Sections~\ref{ssec:linear-structure} and \ref{ssec:cap}.}
  \label{fig:maps-2}
\end{figure}

The only component of $\Phi_{QP}$ involved in the definition of $\phi$
is $Q^{m+1} \to Q^m \otimes P^1$, so $\phi$ reduces level by 1.  The
maps and spaces defined thus far can be organized into the diagram in
\fullref{fig:maps-2}, an expansion of \fullref{fig:maps}.

\begin{lem} \label{lem:phi-chain}\
  \begin{enumerate}
  \item The evaluation map descends to homology and is independent of
    the choice of the representative of \fc.
  \item The cap product map $\phi$ is a chain map of degree $1$.
  \end{enumerate}
\end{lem}

\begin{proof}
  The first part of the lemma stems from the following computation, using
  \fullref{lem:hom-dual}:
  \begin{equation} \label{eqn:d-ev}
    \begin{aligned} 
      \eval \circ \fdf_P(p) &= \langle \fdf_P(p),\fc \rangle \\
      &= \langle p, \fdf_Q (\fc) \rangle  \\
      &= 0.
    \end{aligned}
  \end{equation}
  
  The fact that $\phi$ is a chain map follows from the fact that it is
  the composition of two chain maps (the evaluation map is a chain map
  into the trivial complex $(\zz/2,0)$).  The fact that $\phi$ is a
  degree $1$ map follows from the facts that $\Phi_{QP}$ has degree
  $-1$ (it is derived from $\fdf$) and $\eval$ has degree $2$ (it is
  nonzero only on degree $-2$ generators in $P^1$ dual to the degree
  $1$ class \fc).
\end{proof}

The Morse theoretic motivation suggests that $\phi_*$ inverts
$\eta_*$.  This is indeed the case, up to an application of $\tau_*$.

\begin{prop} \label{prop:phi-inv}
  The inverse of $\eta_*\co  H_*(Q^1, \fdf_Q) \to H_*(P^1, \fdf_P)$ is
  $\tau_* \phi_* = \phi_* \tau_*$.
\end{prop}

\proof
  It suffices to find chain homotopies $H\co  P^1 \to P^1$ and $K\co  Q^1
  \to Q^1$ such that
  \begin{equation} \label{eqn:H}
    \eta \phi \tau + \iota_P = H \fdf_P + \fdf_P H
  \end{equation}
  and 
  \begin{equation} \label{eqn:K} 
    \phi \eta \tau + \iota_Q = K \fdf_Q + \fdf_Q K,
  \end{equation}
  where $\iota_P\co  P^1 \to P^1$ and $\iota_Q\co Q^1 \to Q^1$ are the
  identity maps.

  First, define $H$ by:
  \begin{equation*}
    H = (1 \otimes \eval) \Phi_{PP} \tau.
  \end{equation*}
  Consider the $P^2 \to P^1 \otimes P^1$ component of $\fdf^2 = 0$:
  \begin{equation*}
    \Phi_{PP} \fdf_P + (\fdf_P \otimes 1 + 1 \otimes \fdf_P) \Phi_{PP} +
    (\eta \otimes 1) \Phi_{QP} + 
    (1 \otimes \eta) \Phi_{PQ} = 0.
  \end{equation*}
  Precomposing with $\tau$, post-composing with $(1 \otimes \eval)$,
  and using Lemmas~\ref{lem:translation} and \ref{lem:phi-chain}
  yields
  \begin{equation*}
    H \fdf_P + \fdf_P H + \eta \phi \tau + (1 \otimes \eval \eta)
    \Phi_{PQ} \tau = 0. 
  \end{equation*}
  Let $\iota_P = (1 \otimes \eval \eta) \Phi_{PQ} \tau$; this is
  clearly a chain map since it is a composition of chain
  maps. Further, since $\eta$ acts trivially on $Q^0$, only the
  component of $\Phi_{PQ}$ with image in $P^1 \otimes Q^1$ contributes
  to the definition of $\iota_P$. It now suffices to prove that
  $\iota_P$ is the identity map.

  Let $p \in P^2$ represent a homology class.  Thick disks contribute
  terms of the form $p' \otimes q$ to $\Phi_{PQ} \tau (p)$, where, by
  the K\"unneth theorem, $p'$ represents a homology class in $H_*(P^1,
  \fdf_P)$ and $q$ represents a homology class in $H_*(Q^1, \fdf_Q)$.
  In fact, since a $d$ generator cannot appear as the negative corner
  of a thick disk, $q$ must represent a class in $
  H_*(\overline{Q}^1,\fdf_{\overline{Q}})$.  It follows from
  \fullref{thm:fc-unique} that $\eval \eta (q) = 0$, and hence
  that thick disks not contribute to $\iota_P$.  On the other hand,
  there is always exactly one thin disk that contributes a term of the
  form $p \otimes d$ to $\Phi_{PQ} \tau (p)$. Since
  \fullref{thm:fc-unique} implies that $\eval \eta (d) = 1$, it
  follows that $\iota_P (p) = p$.

  The definition of $K$ and the proof of \eqref{eqn:K} are almost
  identical to the above.  To define $K$, let $\Psi_{QP}\co  Q^2 \to Q^1
  \otimes P^1$ be a component of $\fdf_2$.  Then let
  \begin{equation} \label{eqn:K-def}
    K =  (1 \otimes \eval) \Psi_{QP} \tau.\rlap{\hglue 1.6in \qed}
  \end{equation}

\begin{exa}
  To illustrate the theorem, return to the original observation of
  duality for the knot in \fullref{fig:numbered-5-2}.  As noted at
  the end of \fullref{ssec:linear-structure}, $\eta_*[q_6^1] =
  [p_7^1]$, the latter of which is dual to $[q_7^0]$ with respect to
  the pairing $\langle,\rangle$.  How does $\phi_*$ behave here?

  A computation similar to that in \fullref{ex:5-2-expanded} shows
  that
  \begin{equation*}
    \fdf_2 p_7^2 =  q_6^1 p_9^1 + \text{terms not in } Q^1 \otimes P^1.
  \end{equation*}
  Since $\fc = [q_8+q_9]$, the first term contributes non-trivially
  to $\phi$:
  \begin{equation} \label{eqn:5-2-hom-phi}
    \phi_* [p_7^1] = [q_6^1].
  \end{equation}
  This illustrates that the cap product map $\phi_*$ is the inverse of
  the duality map $\eta_*$.
\end{exa}

\bibliographystyle{gtart}
\bibliography{link}

\end{document}